\documentclass[reqno]{amsart}
\usepackage{amsfonts,amsmath,amssymb,color}

\numberwithin{equation}{section}

\usepackage[babel=true,kerning=true]{microtype}

\usepackage{amsthm,leqno}
\usepackage{hyperref}

\usepackage{amsfonts,amsmath,amssymb,enumerate}
\usepackage{graphicx}
\usepackage{tikz}
\usetikzlibrary{decorations.markings}
\usetikzlibrary{plotmarks}
\usetikzlibrary{patterns}

\newtheorem{theo}{Theorem}
\newtheorem{cl}{Claim}

\newcommand{\eps}{\epsilon}

\begin{document}

\title[Harmonic maps with free boundary]{Regularity and quantification for harmonic maps with free boundary}
\author{Paul Laurain}
\address{Paul Laurain, Department of Mathematics, Building 380, Stanford, California 94305, U.S.A.}
\address{IMJ- Paris 7 
Institut de Math\'ematiques de Jussieu
B\^atiment Sophie Germain
Case 7012
75205 PARIS Cedex 13
France}
\email{laurainp@math.jussieu.fr}

\author{Romain Petrides} 
\address{Romain Petrides, Universit\'e de Lyon, CNRS UMR 5208, Universit\'e Lyon 1, Institut Camille Jordan, 43 bd du 11 novembre 1918, F-69622 Villeurbanne cedex, France.}
\email{romain.petrides@univ-lyon1.fr}

\begin{abstract} 
We prove a quantification result for harmonic maps with free boundary from arbitrary Riemannian surfaces into the unit ball of ${\mathbb R}^{n+1}$ with bounded energy. This generalizes results obtained by Da Lio \cite{DAL2012} on the disc.
\end{abstract}

\maketitle

Let $(M,g)$ be a smooth Riemannian surface with a smooth nonempty boundary with $s$ connected components. We fix $n\ge 2$ and let ${\mathbb B}^{n+1}$ be the unit ball of ${\mathbb R}^{n+1}$. A map $u:\left(M,g\right) \rightarrow \mathbb{B}^{n+1}$ is a smooth harmonic map with free boundary if it is harmonic and smooth up to the boundary, $u(\partial M) \subset \mathbb{S}^n$ and $\partial_{\nu} u$ is parallel to $u$, (or $\partial_{\nu} u \in (T_{u} \mathbb{S}^{n})^\bot$). The energy of such a map is defined as
$$ E(u) = \int_M \left\vert \nabla u\right\vert_g^2\, dv_g = \int_{\partial M} u\cdot\partial_{\nu} u d\sigma_g \hskip.1cm.$$ 
Harmonic maps with free boundary are a natural counterpart of $\mathbb{S}^n$-valued harmonic maps which are critical points of the energy $E$ with the constraint that $\left|u \right|^2 = 1$ on the surface. The difference is that one requires $\left|u \right|^2 = 1$ only on the boundary. 

Harmonic maps with free boundary naturally appear for instance as critical points of Steklov eigenvalues when the metrics stay in a fixed conformal class, as noticed by Fraser and Schoen \cite{FRS2013}. This is the counterpart of $\mathbb{S}^n$-valued harmonic maps which are critical points for the Laplace eigenvalue when the variation of metrics lie in a fixed conformal class (see for instance Petrides \cite{PET2013}). There are also of course strong links between critical points of Steklov eigenvalues, harmonic maps with free boundary and minimal surfaces in a $(n+1)$-ball with free boundary conditions. We refer once again to Fraser-Schoen \cite{FRS2013}.

We prove the following quantification result for harmonic maps with free boundary~:

\begin{theo} \label{t2} Let $u_m : (M,g) \rightarrow \mathbb{B}^{n+1}$ a sequence of harmonic maps with free boundary, i.e $u_m(\partial M) \subset \mathbb{S}^n$ and $u_m$ is parallel to $\partial_{\nu}u_m$ such that
$$ \limsup_{m\to+\infty} \int_{M} \left|\nabla u_m\right|_g^2dv_g < +\infty \hskip.1cm.$$
Then, there is a harmonic map with free boundary $u_{\infty}: M \rightarrow \mathbb{B}^{n+1}$ and
\begin{itemize}
\item $\omega^1,\cdots,\omega^l$ a family of $1/2$-harmonic maps $\mathbb{R} \rightarrow \mathbb{S}^n$,
\item $a_m^1,\cdots,a_m^l$ a family of converging sequences of points on $\partial M$,
\item $\lambda_m^1,\cdots,\lambda_m^l$ a family of sequences of positive numbers all converging to $0$,
\end{itemize}
such that up to the extraction of a subsequence,
$$ u_m \to u_{\infty} \hbox{ in } \mathcal{C}^{\infty}_{loc}( M\setminus \{a_{\infty}^1,\cdots,a_{\infty}^l\} ) \hskip.1cm,$$
and
$$ \int_{\partial M} R_m . \partial_{\nu} R_m \to 0 \hskip.1cm, $$
where we identify $\partial M$ to $s$ copies of $\mathbb{S}^1 = \mathbb{R}\cup \{\infty\}$ :
$$ \partial M = \bigcup_{j=1}^s C_j \hskip.1cm,$$
with $a_{\infty}^i \in C_j \setminus \{\infty\} $ for some $j$,
$$ R_m = u_m - u_{\infty} - \sum_{i=1}^l \omega^i\left(\frac{. - a_m^i}{\lambda_m^i}\right) \hskip.1cm.$$
\end{theo}

In particular, in the space of measures on the boundary, we have the convergence 
$$ u_m. \partial_{\nu} u_m d\sigma_g \rightharpoonup_{\star} u_{\infty}. \partial_{\nu} u_{\infty} d\sigma_g + \sum_{i=1}^l e_i \delta_{a^i_{\infty}} \hskip.1cm,$$
where $e_i$ denotes the energy of the harmonic extension of $\omega^i$ on $\mathbb{R}^2_+$, noticing that a map $\omega$ on the real line $\mathbb{R}\times \{0\}$ is $\frac{1}{2}$-harmonic if and only if his harmonic extension $\omega : \mathbb{R}^2_+ \rightarrow \mathbb{B}^{n+1}$ is harmonic with free boundary. Notice also that a harmonic map with free boundary and finite energy $\omega : \mathbb{R}^2_+ \rightarrow \mathbb{B}^{n+1}$ corresponds to a harmonic map with free boundary on the disc by Claim \ref{remsing} below and the remark which follows. Therefore, $\omega_i$ is automatically conformal and by Fraser and Schoen \cite{FRA2014}, we get that $\omega_i(\mathbb{D})$ is an equatorial plane disc and that the energy of such a map satisfies
$$ e_i = E(\omega_i) = \int_{\mathbb{R}\times \{0\}} \omega_i . (-\partial_t \omega_i)ds \in 2\pi \mathbb{N} \hskip.1cm.$$ 
Thus our result is indeed a quantification result, analogous to those of Sacks-Uhlenbeck \cite{SAC1981}, Parker \cite{PAR1996} for $\mathbb{S}^n$-valued harmonic maps, or of Laurain-Rivi\`ere \cite{LAU2014} for similar equations.

\medskip In the case of the disc, $\left(M,g\right)=\left({\mathbb D},\xi\right)$, this theorem has already been proved by Da Lio \cite{DAL2012}. In this case, the correspondance between harmonic maps with free boundary on the disc with harmonic extensions of $\frac{1}{2}$-harmonic map on the real line is used. Denoting her energy as
$$ E(u) = \int_{\mathbb{R}} \left|\Delta^{\frac{1}{4}}u\right|^2 dx \hskip.1cm,$$
she made use of another variational problem on the real line. This correspondance is no more true for general surfaces.

\medskip The proof of our theorem relies classically on a $\epsilon$-regularity property. Proving it permits also to prove a regularity result on weak harmonic maps with free boundary. A map $u : (M,g) \rightarrow \mathbb{B}^{n+1}$ is called weakly harmonic with free boundary if $ u(x) \in \mathbb{S}^n $ for a.e $x\in\partial M$
and if for any $v\in L^{\infty}\cap H^1(M,\mathbb{R}^{n+1})$ with $v(x) \in T_{u(x)}\mathbb{S}^n$ for a.e $x\in \partial M$,
$$ \int_M  \left\langle\nabla u, \nabla v \right\rangle dv_g = 0\hskip.1cm. $$
Such a map is a critical point for the above energy $E$ with respect to the variations $u_t = \pi(u+tv) $ for any $v\in L^{\infty}\cap H^1(M,\mathbb{R}^{n+1})$ with $v(x) \in T_{u(x)}\mathbb{S}^n$ for a.e $x\in \partial M$, where for $z \in \mathbb{R}^{n+1}$, $\pi(z)$ denotes the the nearest point retraction to $z$ in $\mathbb{B}^{n+1}$. Then we have the following~:

\begin{theo}[Scheven \cite{SCH2006}]\label{t1}
A weakly harmonic map with a free boundary $u : (M,g) \rightarrow \mathbb{B}^{n+1}$ is always smooth until the boundary and thus is a classical harmonic map with free boundary.
\end{theo}

This regularity theorem was originally proved by Scheven \cite{SCH2006}, even in a more general context. This theorem is the analog of that of H\'elein \cite{HEL1996} for $\mathbb{S}^n$-valued harmonic maps.

Another proof, in the case of the disc, was also given by Da Lio-Rivi\`ere \cite{DAL2011}, using the correspondence already explained above. Our proof of this result is an adaptation of the proof of Scheven \cite{SCH2006}. However, we are more careful and precise in the regularity estimate, thus transforming this regularity estimate into an a priori estimate, passing from a regularity result to an $\eps$-regularity result. This is crucial to prove Theorem \ref{t2}.

\medskip As already said, harmonic maps with free boundary in the unit ball naturally appear as critical points of Steklov eigenvalues. Both our theorems are crucial in order to prove existence of regular metrics which maximize $k$-th Steklov eigenvalue on a surface, either with a conformal class constraint or not, as was stressed by Fraser and Schoen \cite{FRS2013}. This is achieved in Petrides \cite{PET2015}.

\medskip Our paper is organized as follows~: section \ref{section-regularity} is devoted to the proof of Theorem \ref{t1}, thanks to an $\epsilon$-regularity result, see Claim 3. We also obtain a crucial result of singularity removability in Claim 4. Our proof is based on the rewriting of the equation with a suitable structure, see Claim 2, permitting to use Wente's inequality, as studied carefully by Rivi\`ere \cite{RIV2012}. Section \ref{section-quantification} is then devoted to the proof of Theorem \ref{t2}. Capitalizing on the $\eps$-regularity result, we are able to prove a no-neck energy result after having described the concentration phenomenon which may appear.

\medskip {\bf Acknowledgements} : The first author was visiting  the Department of Mathematics of Stanford University  when this article was written, he would like to thank it for its hospitality and the excellent working conditions.

\section{Regularity results}\label{section-regularity}

We denote by $\mathbb{D}_r(x)$ the Euclidean disc centered at $x \in \mathbb{R}^2$ of radius $r$ and we let $\mathbb{D}_r = \mathbb{D}_r(0)$ and $\mathbb{D} = \mathbb{D}_1$. On the upper half space, we let $\mathbb{D}_r^+ = \mathbb{D}_r \cap \left(\mathbb{R} \times \mathbb{R}_+ \right)$ and $\mathbb{D}^+ = \mathbb{D}_1^+$. 

We first recall a lemma proved by Scheven \cite{SCH2006}, lemma 3.1, which states that a weakly harmonic map with free boundary with small energy cannot vanish close to the boundary $\partial M$.

\begin{cl}[\cite{SCH2006}] \label{Scheven} There exists $\eps_0$ and $C>0$ such that for any $0<\eps<\eps_0$,  and any weakly harmonic map $u \in H^1(\mathbb{D}^+,\mathbb{B}^{n+1})$ such that $u(\mathbb{R}\times \{0\} ) \subset \mathbb{S}^n$, if
$$ \int_{\mathbb{D}^+} \left|\nabla u\right|^2 \leq \eps \hskip.1cm,$$
then for any $x \in \mathbb{D}_{\frac{1}{2}}^+$,
$$ d(u(x),\mathbb{S}^n) \leq C \eps^{\frac{1}{2}} \hskip.1cm.$$
\end{cl}

To prove the regularity result, we will extend the weakly harmonic map with free boundary $u : (M,g) \rightarrow \mathbb{B}^{n+1}$ by a symmetrization with respect to the boundary. 

Since $u$ is harmonic in the interior of $M$, it is smooth. It remains to prove that $u$ is smooth at the neighbourhood of each point $a \in \partial M$. We take a conformal chart $\phi : U \rightarrow \mathbb{R}^2_+$ centered at $a \in U$ such that
\begin{itemize}
\item $\phi(a) = 0$,
\item $\phi^{\star}(e^{2w}\xi) = g$, where $w$ is some smooth function and $\xi$ is the Euclidean metric,
\item $\phi(\partial M \cap U) = I = (-1,1)\times 0$.
\end{itemize}
We let $\bar{u} = u\circ\phi^{-1}$. Then, by conformal invariance, we have that for any $v \in \mathcal{C}^{\infty}_c(\mathbb{D}^+,\mathbb{R}^{n+1})$, $v(x) \in T_{\bar{u}(x)} \mathbb{S}^n$ for a.e $x\in(-1,1)\times \{0\}$,
$$ \int_{\mathbb{D}^+} \left\langle \nabla v , \nabla \bar{u}\right\rangle = 0 \hskip.1cm.$$
Applying Claim \ref{Scheven}, we can assume up to a dilation  that 
\begin{equation} \label{sch} \forall x\in \mathbb{D}^+ , \left| \bar{u}(x) \right| \geq \frac{1}{2} \hskip.1cm. \end{equation}
In particular, $\bar{u}$ does not vanish and we can define its extension in $\mathbb{D}$ as follows~:
\begin{equation} \label{defutilde} \tilde{u} = \left\{ \begin{array}{ll}
\bar u & \hbox{ on } \mathbb{D}^+ \\
\sigma \circ \bar{u} \circ r & \hbox{ on }  \mathbb{D}^- \end{array} \right. \end{equation}
where for $z\in \mathbb{R}^{n+1}$, $\sigma(z) = \frac{z}{\left|z \right|^2}$ is the inversion w.r.t. the unit sphere $\mathbb{S}^n$, and for $x=(s,t) \in \mathbb{R}^2$, $r(x) = (s,-t)$ and $\mathbb{D}^- = r(\mathbb{D}^+)$.

\begin{cl} \label{weakeq} We have that $\tilde{u} \in H^1(\mathbb{D})$ and satisfies in a weak sense~: for $1\leq j\leq n+1$,
\begin{equation} \label{symeq} - div(A \nabla \tilde{u}_j) = \sum_{i=1}^{n+1} \left\langle X_{i,j} , \nabla \tilde{u}_i \right\rangle \hbox{ in }{\mathbb D} \end{equation}
where $A \in H^{1}(\mathbb{D})$ is defined by
$$ A = \left\{ \begin{array}{ll} 
1 & \mathbb{D}^{+} \\
\frac{1}{\left|\tilde{u}\right|^2} & \mathbb{D}^{-} \end{array} \right. \hskip.1cm, $$
and for $1 \leq i,j \leq n+1$, $X_{i,j} \in L^2(\mathbb{D})$ is defined by
$$ X_{i,j} =  \left\{ \begin{array}{ll} 
0 & \mathbb{D}^{+} \\
2\frac{\tilde{u}_j \nabla \tilde{u}_i - \tilde{u}_i \nabla \tilde{u}_j}{\left|\tilde{u}\right|^4} & \mathbb{D}^{-} \end{array} \right. \hskip.1cm, $$
and satisfies in a weak sense
\begin{equation}\label{divX} div(X_{i,j}) = 0 \hbox{ in }{\mathbb D}\hskip.1cm. \end{equation}
\end{cl}

\medskip {\bf Proof.}
We immediately have that $\tilde{u}$ and $X_{i,j}$ satisfy 
$$ \Delta \tilde{u} = 0 \hbox{ and } div(X_{i,j}) = 0 \hbox{ in } \mathbb{D}^+\setminus I $$
in a classical sense. Direct computations give that 
$$ -div\left( \frac{\nabla \tilde{u}}{\left| \tilde{u}\right|^4}\right) = 2 \frac{\left|\nabla \tilde{u}\right|^2}{\left|\tilde{u}\right|^6} \tilde{u} \hbox{ on } \mathbb{D}^- \setminus I $$
so that  
$$ div(X_{i,j}) = 0 \hbox{ in } \mathbb{D}^- \setminus I  $$
for $1\leq i,j \leq n+1$ and 
$$ -div \left( \frac{\nabla \tilde{u}_j}{\left| \tilde{u}\right|^2}\right) = \sum_{i} \left\langle X_{i,j}, \nabla \tilde{u}_i \right\rangle \hbox{ in } \mathbb{D}^- \setminus I$$
for $1\leq j\leq n+1$, both in a classical sense. It remains to prove that these equations are still true on $\mathbb{D}$ in a weak sense. Thus, for equation (\ref{symeq}), we need to prove that 
\begin{equation} \label{weakeqtildeu} \forall v \in \mathcal{C}_c^{\infty}(\mathbb{D},\mathbb{R}^{n+1}), \int_{\mathbb{D}^+} \left\langle\nabla{\tilde{u}},\nabla v \right\rangle + \int_{\mathbb{D}^-} \left(\frac{\left\langle\nabla{\tilde{u}},\nabla v \right\rangle}{\left|\tilde{u}\right|^2} - \sum_{i,j} \left\langle X_{i,j}, \nabla \tilde{u}_i \right\rangle v_j\right) = 0 \hskip.1cm. \end{equation}
For that purpose, we first remark that, for any $w \in L^{\infty}\cap H^1(\mathbb{D},\mathbb{R}^{n+1})$, we have by direct computations, using the change of variable by the reflection $r$, that
\begin{equation}\label{chvar}
\int_{\mathbb{D}^+} \left\langle \nabla \tilde{u},\nabla w \right\rangle = \int_{\mathbb{D}^-} \left\langle \nabla( \sigma \circ \tilde{u} ), \nabla \left( w\circ r \right) \right\rangle 
\end{equation}
since $\sigma \circ \tilde{u} = \tilde{u}\circ r$. More lengthy but straightforward computations also lead to 
\begin{equation} \label{eqD-} 
\begin{array}{l}
{\displaystyle \int_{\mathbb{D}^-} \left(\frac{\left\langle\nabla{\tilde{u}},\nabla w \right\rangle}{\left|\tilde{u}\right|^2} - \sum_{i,j} \left\langle X_{i,j}, \nabla \tilde{u}_i \right\rangle w_j\right)}\\
{\displaystyle \hskip1cm = \int_{\mathbb{D}^-} \left\langle \nabla \left( \sigma\circ\tilde{u} \right), \nabla \left(w-2\sum_j \frac{\tilde{u}_j w_j}{\vert \tilde{u}\vert^2} \tilde{u}\right)  \right\rangle}\hskip.1cm.
\end{array}
\end{equation}
Indeed,

\begin{equation} \label{eqD-lht} 
\begin{array}{l}
{\displaystyle \frac{\left\langle\nabla{\tilde{u}},\nabla w \right\rangle}{\left|\tilde{u}\right|^2} - \sum_{i,j} \left\langle X_{i,j}, \nabla \tilde{u}_i \right\rangle w_j }\\
{\displaystyle \hskip1cm = \frac{\left\langle\nabla{\tilde{u}},\nabla w \right\rangle}{\left|\tilde{u}\right|^2} - 2 \sum_j \frac{\tilde{u}_j w_j }{\left\vert \tilde{u} \right\vert^4} \left\vert \nabla \tilde{u} \right\vert^2 + \sum_j \frac{ \left\langle \nabla \tilde{u}_j , \nabla \left\vert \tilde{u} \right\vert^2 \right\rangle w_j}{ \left\vert \tilde{u} \right\vert^4 }}  \hskip.1cm,
\end{array}
 \end{equation}
 
 \begin{equation} \label{eqD-rht1} 
\nabla \left( \sigma\circ\tilde{u} \right) = \nabla \left( \frac{\tilde{u} }{ \left\vert \tilde{u} \right\vert^2 } \right) = \frac{ \nabla \tilde{u} }{ \left\vert \tilde{u} \right\vert^2 }  - \frac{ \nabla{ \left( \left\vert \tilde{u} \right\vert^2 \right) \tilde{u}}}{ \left\vert \tilde{u} \right\vert^4 } 
 \end{equation}
and

\begin{equation} \label{eqD-rht2} 
\begin{array}{l}
{\displaystyle    \nabla \left(w-2\sum_j \frac{\tilde{u}_j w_j}{\vert \tilde{u}\vert^2} \tilde{u}\right)  = \nabla w - 2 \sum_j \frac{\tilde{u}_j w_j }{\left\vert \tilde{u} \right\vert^2}  \nabla \tilde{u}      } \\
{\displaystyle \hskip1cm  - 2 \sum_j \frac{\nabla \tilde{u}_j w_j}{ \left\vert \tilde{u} \right\vert^2 }  \tilde{u} - 2 \sum_j \frac{\nabla  w_j \tilde{u}_j}{ \left\vert \tilde{u} \right\vert^2 }  \tilde{u}  + 2\sum_j \tilde{u}_j w_j \tilde{u} \frac{ \nabla  \left\vert \tilde{u} \right\vert^2 }{ \left\vert \tilde{u} \right\vert^4} } \hskip.1cm.
\end{array}
 \end{equation}
 
Let now $v\in \mathcal{C}_c^{\infty}(\mathbb{D},\mathbb{R}^{n+1})$ and let us set 
$$v_e \circ r= \frac{1}{2}\left(v\circ r+v - 2\left(\tilde{u}\cdot v\right)\frac{\tilde{u}}{\left\vert \tilde{u}\right\vert^2}\right)$$
and 
$$v_a\circ r= \frac{1}{2}\left(v\circ r -v+ 2\left(\tilde{u}\cdot v\right)\frac{\tilde{u}}{\left\vert \tilde{u}\right\vert^2}\right)$$
so that $v_a+v_e=v$. Note that $v_a$ and $v_e$ are in $L^{\infty}\cap H^1(\mathbb{D},\mathbb{R}^{n+1})$. Note also that we have 
$$v_e\circ r = v_e - 2\frac{v_e\cdot \tilde{u}}{\left\vert \tilde{u}\right\vert^2} \tilde{u}$$
and
$$v_a \circ r = -\left(v_a -2\frac{v_a\cdot \tilde{u}}{\left\vert \tilde{u}\right\vert^2} \tilde{u}\right)\hskip.1cm.$$
Then we can write, applying \eqref{chvar} and \eqref{eqD-} with $w=v_a$ and $w=v_e$ and using these last equalities, that 
\begin{eqnarray*} 
&&\int_{\mathbb{D}^+} \left\langle\nabla{\tilde{u}},\nabla v \right\rangle + \int_{\mathbb{D}^-} \left(\frac{\left\langle\nabla{\tilde{u}},\nabla v \right\rangle}{\left|\tilde{u}\right|^2} - \sum_{i,j} \left\langle X_{i,j}, \nabla \tilde{u}_i \right\rangle v_j \right)\\
&&\quad =\int_{\mathbb{D}^+} \left\langle\nabla{\tilde{u}},\nabla v_a \right\rangle + \int_{\mathbb{D}^-} \left(\frac{\left\langle\nabla{\tilde{u}},\nabla v_a \right\rangle}{\left|\tilde{u}\right|^2} - \sum_{i,j} \left\langle X_{i,j}, \nabla \tilde{u}_i \right\rangle \left(v_a\right)_j\right) \\
&&\qquad +\int_{\mathbb{D}^+} \left\langle\nabla{\tilde{u}},\nabla v_e \right\rangle + \int_{\mathbb{D}^-} \left(\frac{\left\langle\nabla{\tilde{u}},\nabla v_e \right\rangle}{\left|\tilde{u}\right|^2} - \sum_{i,j} \left\langle X_{i,j}, \nabla \tilde{u}_i \right\rangle \left(v_e\right)_j\right) \\
&&\quad =\int_{\mathbb{D}^+} \left\langle\nabla{\tilde{u}},\nabla v_a \right\rangle -
\int_{\mathbb{D}^-} \langle \nabla \left(  \sigma\circ\tilde{u}\right), \nabla \left(v_a\circ r\right)\rangle \\
&&\qquad + \int_{\mathbb{D}^+} \left\langle\nabla{\tilde{u}},\nabla v_e \right\rangle +
\int_{\mathbb{D}^-} \langle \nabla \left(  \sigma\circ\tilde{u}\right), \nabla \left(v_e\circ r\right)\rangle\\
&&\quad = 2\int_{\mathbb{D}^+} \left\langle\nabla{\tilde{u}},\nabla v_e \right\rangle\hskip.1cm.
\end{eqnarray*}
Noticing now that if $x\in\mathbb{R}\times \{0\}$,
$$ v_e(x) = v(x) - \left(\tilde{u}(x)\cdot v(x) \right) \tilde{u}(x) \in T_{\tilde{u}(x)} \mathbb{S}^n \hskip.1cm,$$
and recalling that $u$ is weakly harmonic with free boundary, we know that
$$ \int_{\mathbb{D}^+} \left\langle \nabla \tilde{u}, \nabla v_e \right\rangle = 0 $$
which clearly ends the proof of \eqref{weakeqtildeu}.

It remains to prove (\ref{divX}), that is
$$ \forall f\in \mathcal{C}_c^{\infty}(\mathbb{D},\mathbb{R}^{n+1}), \int_{\mathbb{D}} \left\langle \nabla f, X_{i,j} \right\rangle = 0\hskip.1cm. $$
Let $f\in \mathcal{C}_c^{\infty}(\mathbb{D},\mathbb{R})$ and $1\leq i,j \leq n+1$ be such that $i\neq j$. Then
\begin{eqnarray*} 
\frac{1}{2}\int_{\mathbb{D}} \left\langle \nabla f, X_{i,j} \right\rangle & = & \int_{\mathbb{D}^-} \left\langle \frac{\tilde{u}_j \nabla \tilde{u}_i - \tilde{u}_i \nabla \tilde{u}_j}{\left| \tilde{u}\right|^4} , \nabla f \right\rangle \\
& = & \int_{\mathbb{D}^+} \left\langle \frac{(\sigma \circ \tilde{u})_j \nabla (\sigma\circ \tilde{u})_i - (\sigma \circ \tilde{u})_i \nabla (\sigma \circ \tilde{u})_j}{\left| \tilde{u}\circ \sigma \right|^4} , \nabla (f\circ r) \right\rangle \\
& = & \int_{\mathbb{D}^+} \left\langle \tilde{u}_j \nabla \tilde{u}_i - \tilde{u}_i \nabla \tilde{u}_j , \nabla (f\circ r) \right\rangle \\
& = & \int_{\mathbb{D}^+} \left\langle \nabla v_{i,j}, \nabla \tilde{u} \right\rangle \\
& =& 0 \hskip.1cm,\\ 
\end{eqnarray*}
where $v_{i,j} \in \mathcal{C}_c^{\infty}\left( \mathbb{D}, \mathbb{R}^{n+1} \right)$ is defined by
$$ (v_{i,j})_k = f\circ r \times \left\{ \begin{array}{ll}
\tilde{u}_j  & \hbox{ if } k=i \\
- \tilde{u}_i  & \hbox{ if } k=j \\
0 & \hbox{ otherwise }
  \end{array} \right. $$
and $\left\langle v_{i,j}(x),\tilde{u}(x) \right\rangle =0 $ for $x\in \mathbb{D}$. This ends the proof of the claim. 
  
\hfill $\diamondsuit$

\medskip Notice that this construction is similar to Scheven's one \cite{SCH2006}. However, in Claim \ref{weakeq}, we give a suitable form to the equation that satisfies the symmetrized map $\tilde{u}$ : we use its structure to prove an $\epsilon$-regularity result that will be useful for the second part of the paper. This type of equations was intensively studied by Rivi\`ere (see \cite{RIV2012}).

\begin{cl} \label{epsreg} There is $\eps_1>0$ and a constant $C_k$ such that if a weakly harmonic map with free boundary $u$ satisfies at the neighbourhood of $a \in \partial M$
$$ \int_{\mathbb{D}^+} \left| \nabla \tilde{u}\right|^2 \leq \eps_1 \hskip.1cm,$$
then $\tilde{u} \in \mathcal{C}^{\infty}(\mathbb{D}_{\frac{1}{2}}^+)$ and for any $k\geq 0$,
$$ \left\| \nabla\tilde{u} \right\|_{\mathcal{C}^{k}(\mathbb{D}_{\frac{1}{2}}^+)} \leq C_k \left\| \nabla\tilde{u} \right\|_{L^{2}(\mathbb{D}^+)}\hskip.1cm. $$
\end{cl}

\medskip {\bf Proof.}
We fix $0<\eps_2<\eps_0$ that we shall choose later, where $\eps_0$ is given by Claim \ref{Scheven} and we assume that
$$ \int_{\mathbb{D}} \left| \nabla \tilde{u}\right|^2 \leq \eps_2 \hskip.1cm.$$
Since $X_{i,j} \in L^2(\mathbb{D})$ satisfies (\ref{divX}), that is $div(X_{i,j}) = 0$ in a weak sense in $\mathbb{D}$, there is a function $B_{i,j} \in H^1(\mathbb{D})$ such that $X_{i,j} = \nabla^{\perp} B_{i,j}$. Then, $\tilde{u} \in H^{1}(\mathbb{D})$ satisfies the equations
$$ \left\{ \begin{array}{ll}
- div(A\nabla\tilde{u}_j) = \sum_i \left\langle\nabla^\perp B_{i,j}, \nabla \tilde{u}_i \right\rangle \\
rot(A \nabla \tilde{u}_i) = \left\langle \nabla^{\perp} A,\nabla \tilde{u}_i \right\rangle
\end{array}  \right.$$
Let $p\in\mathbb{D}_{\frac{1}{2}}$ and $0<r<\frac{1}{2}$. We let $C \in H^1(\mathbb{D})$ be such that
\begin{equation} \label{eqC} \left\{ \begin{array}{ll} 
\Delta C_j = \sum_i \left\langle \nabla^{\perp}B_{i,j} , \nabla \tilde{u}_i \right\rangle & \hbox{ in } \mathbb{D}_r(p) \\
C = 0 & \hbox{ on } \partial \mathbb{D}_r(p)
\end{array} \right. \end{equation} 
Since $div(A\nabla\tilde{u} - \nabla C) = 0$, there is $D\in H^{1}(\mathbb{D})$ such that $\nabla^{\perp} D = A\nabla\tilde{u} - \nabla C $. We set $D = \phi + v$ where $v$ is harmonic and $\phi$ satisfies
\begin{equation} \label{eqphi} \left\{ \begin{array}{ll} 
\Delta \phi_i = - \left\langle \nabla^{\perp} A , \nabla \tilde{u}_i \right\rangle & \hbox{ in } \mathbb{D}_r(p) \\
\phi = 0 & \hbox{ on } \partial \mathbb{D}_r(p)
\end{array} \right. \end{equation} 
Wente's theorem applied to \eqref{eqC} and \eqref{eqphi} gives the estimates
\begin{equation} \label{wenteC} \left\|\nabla C\right\|_{L^2(\mathbb{D}_r(p))} \leq K_0 \left\|\nabla \tilde{u}\right\|_{L^2(\mathbb{D}_r(p))}^2 \hbox{ and }\end{equation}
\begin{equation} \label{wentephi} \left\|\nabla \phi\right\|_{L^2(\mathbb{D}_r(p))} \leq K_0 \left\|\nabla \tilde{u}\right\|_{L^2(\mathbb{D}_r(p))}^2 \hskip.1cm , \end{equation}
with $K_0$ a universal constant. Here we used \eqref{sch}. Moreover, by Rivi\`ere \cite{RIV2012}, lemma VII.1, 
\begin{equation} \label{harmest} \left\|\nabla v\right\|_{L^2(\mathbb{D}_{\frac{r}{16}}(p))} \leq \frac{1}{16} \left\|\nabla v\right\|_{L^2(\mathbb{D}_{r}(p))} \hskip.1cm.\end{equation}
Then we have by Young's inequalities and \eqref{wenteC}, \eqref{wentephi} that
\begin{eqnarray*}
\left\| \nabla\tilde{u} \right\|_{L^2(\mathbb{D}_{\frac{r}{16}}(p))}^2 & = & \left\| A^{-1} A \nabla\tilde{u} \right\|_{L^2(\mathbb{D}_{\frac{r}{16}}(p))}^2 \\
& \leq & 2 \left\| A^{-1} \right\|_{\infty}^2 \left( \left\| \nabla C \right\|^2_{L^{2}(\mathbb{D}_{\frac{r}{16}}(p))} + \left\| \nabla D \right\|^2_{L^{2}(\mathbb{D}_{\frac{r}{16}}(p))} \right) \\
& \leq & 2 \left\| A^{-1} \right\|_{\infty}^2 \left( 3 K_0^2 \left\| \nabla \tilde{u} \right\|^4_{L^{2}(\mathbb{D}_{r}(p))} + 2 \left\| \nabla v \right\|^2_{L^{2}(\mathbb{D}_{\frac{r}{16}}(p))} \right)\hskip.1cm.
\end{eqnarray*}
And an integration by parts on $\mathbb{D}_r(p)$ gives that 
$$ \left\| A \nabla \tilde{u} \right\|_{L^2(\mathbb{D}_r(p))}^2 = \left\| \nabla C \right\|_{L^2(\mathbb{D}_r(p))}^2+\left\|  \nabla D \right\|_{L^2(\mathbb{D}_r(p))}^2$$
since $C=0$ on $\partial {\mathbb D}_r\left(p\right)$. Now we also have that 
$$\left\|  \nabla v \right\|_{L^2(\mathbb{D}_r(p))}^2 = \left\|  \nabla D \right\|_{L^2(\mathbb{D}_r(p))}^2 - \left\|  \nabla \varphi\right\|_{L^2(\mathbb{D}_r(p))}^2\le \left\|  \nabla D \right\|_{L^2(\mathbb{D}_r(p))}^2 $$
so that, since $A(x)\leq 1$ for a.e $x\in\mathbb{D}$, we can write that 
$$
\left\| \nabla v \right\|_{L^2(\mathbb{D}_{\frac{r}{16}}(p))} \leq  \frac{1}{16} \left\| \nabla v \right\|_{L^2(\mathbb{D}_{r}(p))} 
\leq  \frac{1}{16} \left\| \nabla D \right\|_{L^2(\mathbb{D}_{r}(p))} 
\leq  \frac{1}{16} \left\| \nabla \tilde u \right\|_{L^2(\mathbb{D}_{r}(p))}
$$
and we finally get that 
$$ \left\| \nabla\tilde{u} \right\|_{L^2(\mathbb{D}_{\frac{r}{16}}(p))}^2 \leq 6K_0^2 \left\| A^{-1} \right\|_{\infty}^2 \left\| \nabla\tilde{u} \right\|_{L^2(\mathbb{D}_{r}(p))}^4 + \frac{1}{64} \left\| A^{-1} \right\|_{\infty}^2 \left\| \nabla\tilde{u} \right\|_{L^2(\mathbb{D}_{r}(p))}^2 \hskip.1cm.$$
Since $\left\| A^{-1} \right\|_{\infty}^2 \le 16$ thanks to \eqref{sch}, up to choose $\eps_2=\frac{1}{4\times 96 K_0^2}$, we  get that 
$$ \left\| \nabla\tilde{u} \right\|_{L^2(\mathbb{D}_{\frac{r}{16}}(p))}^2 \leq \frac{1}{2} \left\| \nabla\tilde{u} \right\|_{L^2(\mathbb{D}_{r}(p))}^2 $$
for any $p\in \mathbb{D}_{\frac{1}{2}}$ and $0<r<\frac{1}{2}$. Thanks to Morrey estimates, see page 50 of \cite{RIV2012}, and the elliptic regularity on the equation, knowing that $\left| \nabla^{\perp} B\right|^2  \le K \left| \nabla \tilde{u}\right|^2$ almost everywhere for some constant $K$, we get a constant $C$ independent of $\tilde{u}$ such that
\begin{equation} \label{holderest} \left\| \nabla\tilde{u} \right\|_{\mathcal{C}^{1,\gamma}\left(\mathbb{D}_{\frac{1}{2}}\right)} \leq C \left\| \nabla\tilde{u} \right\|_{L^{2}(\mathbb{D})}\hskip.1cm. \end{equation}
Since by (\ref{defutilde}),
$$ \int_{\mathbb{D}^-} \left|\nabla \tilde{u}\right|^2 = \int_{\mathbb{D}^+} \frac{\left|\nabla \tilde{u}\right|^2}{\left|\tilde{u}\right|^4} \hskip.1cm, $$
using (\ref{sch}), and setting $\eps_1 = \frac{\eps_2}{17}$, we get
$$ \int_{\mathbb{D}^+}  \left|\nabla \tilde{u}\right|^2 \leq \eps_1 \Rightarrow \int_{\mathbb{D}}  \left|\nabla \tilde{u}\right|^2 \leq \eps_2 \hskip.1cm. $$
Finally, by elliptic regularity theory, we can bootstrap (\ref{holderest}) thanks to the equation on the half space,
$$ 
\left\{ \begin{array}{ll}
\Delta \tilde{u} = 0 & \hbox{ on } \mathbb{D^+} \\
-\partial_{t} \tilde{u} = \left( \tilde{u} . (-\partial_{t} \tilde{u}) \right)  \tilde{u} & \hbox{ on } I \end{array} \right.
$$
and get the claim.

\hfill $\diamondsuit$

Theorem \ref{t1} of course follows from this claim.

\medskip Thanks to Claim \ref{epsreg}, we also have a result of removability of  singularities for harmonic maps with free boundary which will be useful in the next section.

\begin{cl} \label{remsing} Let $u: M \setminus \{a\} \rightarrow \mathbb{B}^{n+1}$ with finite energy be  such that for any $v \in H^1(M,\mathbb{R}^{n+1})\cap L^{\infty}$, $supp (v) \subset M \setminus \{a\}$ and $v(x) \in T_{u(x)}\mathbb{S}^n$ for a.e $x \in \partial M$,
$$ \int_M \left\langle \nabla u,\nabla v \right\rangle_g dv_g = 0 \hskip.1cm.$$
Then $u$ extends to a harmonic map with free boundary $\bar{u} : M \rightarrow \mathbb{B}^{n+1}$.
\end{cl}

\medskip {\bf Proof.} First, using Claim \ref{epsreg}, it is clear that $u$ is smooth outside of $a$. We use in the sequel the same notation as above since the problem is purely local. Thus we can assume that $M = {\mathbb D}^+$, that the metric is Euclidean and that $a=0$. By a direct scaling argument, using again Claim \ref{epsreg} and standard elliptic regularity theory for harmonic maps in the inside of ${\mathbb D}^+$, there exists a constant $C>0$ such that 
$$\sup_{{\mathbb D}_{\frac{\vert p\vert}{4}}\left(p\right)\cap {\mathbb D}^+} \left\vert x\right\vert^2 \left\vert \nabla u\right\vert^2 \le C \int_{{\mathbb D}_{\frac{\vert p\vert}{2}}\left(p\right)\cap {\mathbb D}^+} \left\vert \nabla u\right\vert^2$$
for all $p\in {\mathbb D}^+$ with $\vert p\vert \le \frac{1}{2}$ as soon as 
$$\int_{{\mathbb D}_{\frac{\vert p\vert}{2}}\left(p\right)\cap {\mathbb D}^+} \left\vert \nabla u\right\vert^2 \le \eps_1\hskip.1cm.$$
Since $\nabla u\in L^2\left({\mathbb D}^+\right)$ by assumption, we deduce that 
\begin{equation}\label{estfaible}
\sup_{{\mathbb D}_r^+} \vert x\vert \left\vert \nabla u \right\vert \to 0\hbox{ as }r\to 0\hskip.1cm.
\end{equation}
Let $v\in C_c^\infty\left({\mathbb D}^+\right)$ be such that for all $x\in\mathbb{R}\times \{0\}$, 
$v(x)\in T_{u(x)} \mathbb{S}^n$. Then we have, integrating by parts, that 
$$\int_{{\mathbb D}^+\setminus {\mathbb D}^+_r} \langle \nabla u,\nabla v\rangle = \int_{\partial {\mathbb D}^+_r} - \partial_\nu u \cdot v\hskip.1cm.$$
Using \eqref{estfaible} and the fact that $\nabla u \in L^2\left({\mathbb D}^+\right)$, we can pass to the limit as $r\to 0$ to obtain that $u$ is in fact a weak harmonic map with free boundary. It is thus regular thanks to Theorem \ref{t1} we just proved. 

\hfill $\diamondsuit$

\medskip Notice that thanks to Claim \ref{remsing}, we have a correspondance between harmonic maps  with free boundary $u : \mathbb{D} \rightarrow \mathbb{B}^{n+1}$ and $v : \mathbb{R}_+^2 \rightarrow \mathbb{B}^{n+1}$, thanks to $f : \mathbb{R}^2_+ \rightarrow \mathbb{D}\setminus \{(0,1)\} $, the conformal map defined by $f(z) = \frac{z-i}{z+i} $. 

Finally, Claim \ref{epsreg} reveals an energy gap for harmonic maps with free boundary on discs: if a harmonic map with free boundary $\omega: \mathbb{R}^2_+ \rightarrow \mathbb{B}^{n+1}$ satisfies
$$ \int_{\mathbb{R}^2_+} \left|\nabla \omega\right|^2 \leq \eps_1 $$
then $\omega$ is a constant map. Indeed, by Claim \ref{epsreg} and an obvious scaling argument, we get that 
$$ \left\| \nabla\omega \right\|_{\mathcal{C}^0(\mathbb{D}_R^+)} \leq \frac{C_0}{R} \left\|\nabla \omega\right\|_{L^2} $$
for all $R>0$ for some fixed constant $C_0$. Letting $R$ go to $+\infty$ gives that $\omega$ is constant.

\section{The quantification phenomenon}\label{section-quantification}

We aim at proving Theorem \ref{t2}. 

\medskip {\bf Step 1 : Points of concentration.}

Since the energy of the sequence $(u_m)$ of harmonic maps with free boundary is bounded, we only have a finite number of points, denoted by $a^1,\cdots,a^q$ such that
\begin{equation}\label{concentpts} \forall r >0, \limsup_{m\to +\infty} \int_{B_g(a^i,r)} \left|\nabla u_m \right|_g^2 dv_g > \eps_1 \hskip.1cm. \end{equation}

Notice that $a^i \in \partial M$ for any $1\leq i \leq q$. Indeed, if $a^i \in M \setminus \partial M$, then, since $u_m$ is harmonic, elliptic regularity theory gives a constant $C$ independent of $m$ such that
$$ \left\| \nabla u_m \right\|_{\mathcal{C}^0(B_g(a^i,\frac{\delta}{2}))} \leq C \left\|  u_m \right\|_{W^{1,2}} \hskip.1cm,$$
where $\delta = d(a^i,\partial M)>0$, so that $(\nabla u_m)$ is uniformly bounded on $B_g(a^i,\frac{\delta}{2})$. This contradicts (\ref{concentpts}).

By $\eps$-regularity around each point of $\partial M\setminus \{a^1,\cdots,a^q \}$, (see Claim \ref{epsreg}), we get that
\begin{equation} \label{strongconv} u_m \to u_{\infty} \hbox{ in } \mathcal{C}_{loc}^1(M \setminus \{a^1,\cdots, a^q \}) \hbox{ as } m\to +\infty \hskip.1cm, \end{equation}
where $u_{\infty}$ satisfies the hypothesis of Claim \ref{remsing} so that $u_{\infty}$ extends to an harmonic map with free boundary in $\mathcal{C}^{\infty}(M,\mathbb{R}^{n+1})$.

\medskip {\bf Step 2 : Blow-up around $a^i \in \partial M$.}

We take a conformal chart $\phi_i : U_i \rightarrow \mathbb{R}^2_+$ centered at $a^i \in U_i$ such that
\begin{itemize}
\item $\phi_i(a^i) = 0$,
\item $\phi_i^{\star}(e^{2w_i}\xi) = g$, where $w_i$ is some smooth function and $\xi$ is the Euclidean metric,
\item $\phi_i(\partial M \cap U_i) \subset \mathbb{R}\times \{0\}$.
\end{itemize}

We fix $1\leq i \leq q$ and we let $\bar{u}^i_m = u_m \circ \phi_i^{-1}$ and $\bar{u}_{\infty}^i = u_{\infty} \circ \phi_i^{-1}$. We choose $r_i>0$ small enough such that for any $j\neq i$, $\phi_i(a^j) \notin \mathbb{D}^+_{r_i}$ and at the neighbourhood of $a^i \in \partial M$, 
\begin{equation} \label{chooser} \int_{\mathbb{D}^+_{r_i}} \left|\nabla \bar{u}^i_{\infty}\right|^2 < \frac{\eps_1}{4}\hskip.1cm. \end{equation}
Since $a^i$ satisfies (\ref{concentpts}), we can take $\lambda_m^i$ such that
\begin{equation} \label{chooselambda} \int_{\mathbb{D}^+_{r_i} \setminus \mathbb{D}^+_{\lambda_m^i}} \left|\nabla \bar{u}^i_m  \right|^2 = \frac{\eps_1}{2} \hskip.1cm. \end{equation}
Notice that 
\begin{equation} \label{lambdato0} \lambda_m^i \to 0 \hbox{ as } m\to +\infty \hskip.1cm. \end{equation}
Indeed, if $\displaystyle \limsup_{m\to +\infty} \lambda_m^i >0$, by the definition of $r_i$, and (\ref{strongconv}), passing to the limit in (\ref{chooselambda}) would contradict (\ref{chooser}).

We set for $x\in \mathbb{R}^2_+$
$$ \tilde{u}_m^i(x) = \bar{u}_m^i ( \lambda_m^i x) \hskip.1cm.$$

Since $\tilde{u}_m^i$ is harmonic with finite energy, elliptic estimates prove that, $\tilde{u}_m^i$ does not concentrate on $\mathbb{R}\times(0,+\infty)$. Moreover, by (\ref{chooselambda}) and Claim \ref{epsreg}, $\tilde{u}_m^i$ does not concentrate on $\mathbb{R}^2_+ \setminus (-1,1)\times \{0\}$. Therefore, outside some concentration points $a^{i,1},\cdots,a^{i,q_i} \in (-1,1)\times \{0\}$, we have
\begin{equation} \label{convutilde} \tilde{u}_m^i \to \tilde{u}_{\infty}^i \hbox{ in } \mathcal{C}_{loc}^1(\mathbb{R}^2_+\setminus \{a^{i,1},\cdots,a^{i,q_i}\}) \hbox{ as } m\to+\infty \hskip.1cm. \end{equation}
Let $f : \mathbb{R}^2_+ \rightarrow \mathbb{D} $ the conformal map defined by $f(z) = \frac{z-i}{z+i} $. Then $\tilde{u}^i_{\infty} \circ f^{-1}$ satisfies the hypotheses of Claim \ref{remsing} on $\mathbb{D} \setminus \{1,f(a^{i,1}),\cdots,f(a^{i,q_i})\}$ so that $\tilde{u}^i_{\infty} \circ f^{-1}$ extends to an harmonic map with free boundary in $\mathcal{C}^{\infty}(\mathbb{D},\mathbb{R}^{n+1})$.
Thanks to (\ref{chooselambda}), (\ref{lambdato0}) and Claim \ref{noneckenergy}, we have that 
$$ \lim_{R\to+\infty} \lim_{m\to+\infty} \int_{\mathbb{D}^+_{\frac{r_i}{\lambda_m^i R}} \setminus \mathbb{D}^+_{R}} \left|\nabla\tilde{u}_m^i \right|^2 = 0 $$ 
so that by (\ref{chooselambda}) and (\ref{convutilde}) for $R$ large enough,
\begin{equation} \label{energystep2} \int_{\mathbb{D}^+_R \setminus \mathbb{D}^+} \left| \nabla \tilde{u}_{\infty}^i \right|^2 > \frac{\eps_1}{4} \hskip.1cm. \end{equation}
In particular $\tilde{u}_{\infty}^i$ is a non-constant function, which is to a $\frac{1}{2}$-harmonic map on the boundary (one of the $\omega^j$'s given by Theorem \ref{t2}).

\medskip {\bf Step 3 : Iteration.}

As a classical bubble tree extraction (see \cite{PAR1996}), we have two cases: Either there are concentration points and we go back to Step 2 at the neighbourhood of each $a^{i,j}$. Or there is no concentration points for the sequence $(\tilde{u}_m^i)$ (that is $q_i = 0$ in Step 2) and the process stops.

This process has to stop since at every new concentration point, we get a bubble whose energy is at least $\frac{\eps_1}{4}$ by (\ref{energystep2}).

\medskip Finally, we state a no-neck-energy lemma which concludes the proof of Theorem \ref{t2}.

\begin{cl} \label{noneckenergy} Let $(\lambda_m)$ be a sequence of positive numbers converging to $0$. Let $(u_m)$ be a sequence of harmonic maps on $\mathbb{D}^+$ with uniformly bounded energy and free boundary on $(-1,1)\times\{0\}$ such that 
\begin{equation} \label{assumpnoneck} \int_{\mathbb{D}^+ \setminus \mathbb{D}^+_{\lambda_m}} \left|\nabla u_m\right|^2 \leq \frac{\eps_1}{2} \hskip.1cm, \end{equation}
Then,
\begin{equation} \label{limgrad} \lim_{R\to +\infty} \lim_{m\to +\infty} \int_{\mathbb{D}^+_{\frac{1}{R}} \setminus \mathbb{D}^+_{\lambda_m R}} \left|\nabla u_m\right|^2 = 0 \end{equation}
and
\begin{equation} \label{limboundary} \lim_{R\to +\infty} \lim_{m\to +\infty} \int_{(-\frac{1}{R},\frac{1}{R}) \setminus (-\lambda_m R,\lambda_m R)} u_m . \partial_t u_m = 0 \hskip.1cm. \end{equation}
\end{cl}

\medskip {\bf Proof.}

We set ${\mathbb A}_{m,R} = \mathbb{D}^+_\frac{1}{R}\setminus \mathbb{D}^+_{\lambda_m R}$, $I_{m,R} = {\mathbb A}_{m,R} \cap \left(\mathbb{R}\times \{0\}\right)$ and
$$\delta_{m,R} = \max_{z \in {\mathbb A}_{m,R}} \left|z\right| \left|\nabla u_m \right|(z) \hskip.1cm.$$

\medskip {\bf Step 1:}We have that
\begin{equation} \label{limdeltamR} \lim_{R\to +\infty} \limsup_{m\to+\infty}\delta_{m,R} = 0\hskip.1cm.  
\end{equation}

\medskip {\it Proof of Step 1} - Notice that $R\mapsto \delta_{m,R}$ and $R \mapsto \displaystyle{\limsup_{m\to+\infty}\delta_{m,R}}$ are nonincreasing. We proceed by contradiction, assuming \eqref{limdeltamR} is false. Then there exists a subsequence $\left(m_\alpha\right)_{\alpha\ge 1}$ converging to $+\infty$ such that 
\begin{equation}\label{limdeltamR1}
\delta_{m_\alpha,\alpha} \ge \eps_0>0
\end{equation}
for some $\eps_0>0$ fixed. Let $z_\alpha \in {\mathbb A}_{m_\alpha,\alpha}$ be such that $\delta_{m_\alpha,\alpha} = \left|z_\alpha\right| \left|\nabla u_{m_\alpha} \right|\left(z_\alpha\right)$. It is clear that 
$z_\alpha\to 0$ and $\frac{\left|z_\alpha\right|}{\lambda_{m_\alpha}} \to +\infty$ as $\alpha\to +\infty$. We let $\overline{u}_\alpha(x) = u_{m_\alpha}\left(\left|z_\alpha\right| x\right)$ so that, by Claim \ref{epsreg},
$$ \overline{u}_\alpha \to \overline{u}_{\infty} \mbox{ in } \mathcal{C}_{loc}^1(\mathbb{R}^2_+\setminus\{0\}) \hbox{ as } \alpha\to +\infty $$
where $\overline{u}_{\infty}$ is harmonic with free boundary. Then, since, after passing to a subsequence, 
$$ \left|z_\alpha\right| \left|\nabla u_{m_\alpha} \left(z_\alpha\right)\right| = \left|\nabla \overline{u}_\alpha \left(\frac{z_\alpha}{\left|z_\alpha\right|}\right)\right| \to \left|\nabla \overline{u}_{\infty}(z)\right| \hbox{ as } \alpha\to +\infty $$
where $\displaystyle z=\lim_{\alpha\to +\infty}\frac{z_\alpha}{\left|z_\alpha\right|}$, we get thanks to \eqref{limdeltamR1} that $\left|\nabla \overline{u}_{\infty}(z)\right|\ge \eps_0$. By assumption (\ref{assumpnoneck}), $\left\|\nabla \overline{u}_{\infty}\right\|_{L^2(\mathbb{R}^2_+)}^2 \leq \frac{\eps_1}{2}$ so that by the remark at the end of Section 1, $\overline{u}_{\infty}$ should be constant. This is a contradiction which ends the proof of this step.

\medskip {\bf Step 2:}
\begin{equation} \label{steplim} \lim_{R\to +\infty} \limsup_{m\to +\infty} \left\| \nabla u_m \right\|_{L^{2,\infty}\left({\mathbb A}_{m,R}\right)} = 0 \end{equation}
 
\medskip {\it Proof of Step 2} - We easily check that $\left\| \frac{1}{\left|x\right|}\right\|_{L^{2,\infty}\left({\mathbb A}_{m,R}\right)} \leq \sqrt{\pi}$ for any $m$ and $R$, so that
$$ \left\| \nabla u_m \right\|_{L^{2,\infty}\left({\mathbb A}_{m,R}\right)} \leq  \sqrt{\pi} \delta_{m,R} $$
and we get (\ref{steplim}) thanks to Step 1.

\medskip {\bf Step 3:}
\begin{equation} \label{stepangular} \limsup_{R\to +\infty} \limsup_{m\to +\infty} \left\| \nabla_{\theta} u_m \right\|_{L^{2,1}\left({\mathbb A}_{m,R}\right)} < +\infty \end{equation}

\medskip {\it Proof of Step 3} - Here we use the symmetrization process given by Claim \ref{weakeq}. And such estimates on the angular derivative for solutions of this type of equations were obtained in Laurain-Rivi\`ere \cite{LAU2014}. First, we have to ensure that $\left|u_m(x)\right|$ does not vanish for $x\in {\mathbb A}_{m,R}$. For $x \in {\mathbb A}_{m,R}$, we have that
$$ \left\vert u_m(x) - u_m(y)\right\vert \le \frac{\pi}{2} \vert x\vert\sup_{\vert z\vert = \vert x\vert}\left\vert \nabla u_m\right\vert$$
for some $y\in \left(-1,1\right) \times \{0\}$. Since $\vert u(y)\vert=1$, we deduce that 
$$ \left| u_m(x) \right| \geq 1 - \frac{\pi}{2} \delta_{m,R} $$
Since $R \mapsto \delta_{m,R}$ are decreasing for every $m$ and thanks to (\ref{limdeltamR}), we deduce that there exists $R_0$ and $m_0$ such that for all $R\ge R_0$ and all $m\ge m_0$, 
$$ \left| u_m \right| \geq \frac{1}{2} \hbox{ in } {\mathbb A}_{m,R}\hskip.1cm.$$
Up to a scaling and since we are interested only in large $m$ and large $R$ and in order to simplify the notations, we may assume that $R_0=1$ and that $m_0=1$. We set now  
\begin{equation} \label{defutildem} \tilde{u}_m = \left\{ \begin{array}{ll}
u_m & \hbox{ in } \mathbb{D}^+ \setminus \mathbb{D}_{\lambda_m}^+ \\
\sigma \circ u_m \circ r & \hbox{ in }  \mathbb{D}^- \setminus \mathbb{D}_{\lambda_m}^- \end{array} \right. \end{equation}
Applying the computations of Claim \ref{weakeq}, we get that, for $1\leq j\leq n+1$,
\begin{equation} \label{symeqm} - div(A_m \nabla \tilde{u}_m^j) = \sum_{i=1}^{n+1} \left\langle X_m^{i,j} , \nabla \tilde{u}_m^i \right\rangle  \end{equation}
in a weak sense where $A_m$ is defined by
$$ A_m = \left\{ \begin{array}{ll} 
1 & \mathbb{D}^{+} \setminus \mathbb{D}_{\lambda_m}^+ \\
\frac{1}{\left|\tilde{u}_m\right|^2} & \mathbb{D}^{-} \setminus \mathbb{D}_{\lambda_m}^- \end{array} \right. \hskip.1cm, $$
and for $1 \leq i,j \leq n+1$, $X_m^{i,j}$ is defined by
$$ X_m^{i,j} =  \left\{ \begin{array}{ll} 
0 & \mathbb{D}^{+} \setminus \mathbb{D}_{\lambda_m}^+ \\
2\frac{\tilde{u}_m^j \nabla \tilde{u}_m^i - \tilde{u}_m^i \nabla \tilde{u}_m^j}{\left|\tilde{u}_m\right|^4} & \mathbb{D}^{-} \setminus \mathbb{D}_{\lambda_m}^- \end{array} \right. \hskip.1cm, $$
and satisfies in a weak sense
\begin{equation}\label{divXm} div(X_m^{i,j}) = 0 \hskip.1cm. \end{equation}
We also have that for $\lambda_m < r < 1$,
\begin{eqnarray*} 
\frac{1}{2}\int_{\partial\mathbb{D}_r} \left\langle X_m^{i,j},\nu \right\rangle & = & \int_{\partial\mathbb{D}_r^-} \left\langle \frac{\tilde{u}_m^j \nabla \tilde{u}_m^i - \tilde{u}_m^i \nabla \tilde{u}_m^j}{\left| \tilde{u}_m\right|^4} , \nu \right\rangle \\
& = & \int_{\partial \mathbb{D}_r^+} \left\langle \frac{(\sigma \circ u_m)^j \nabla (\sigma\circ u_m)^i - (\sigma \circ u_m)^i \nabla (\sigma \circ u_m)^j}{\left| \sigma\circ u_m\circ \right|^4} , \nu \right\rangle \\
& = & \int_{\partial \mathbb{D}_r^+} \langle u_m^j \partial_{\nu} u_m^i - u_m^i \partial_{\nu} u_m^j, \nu\rangle \\
& = & \int_{[-r,r]\times \{0\}} \left(u_m^j \partial_{t} u_m^i - u_m^i \partial_{t} u_m^j \right) \\
&& \quad- \int_{\mathbb{D}_r^+} div \left(u_m^j \nabla u_m^i - u_m^i \nabla u_m^j \right) \\
& =& 0 \hskip.1cm,\\ 
\end{eqnarray*}
since $\Delta u_m = 0$ and $\partial_t u_m$ is parallel to $u_m$ on $[-r,r]\times \{0\}$. From this and \eqref{divXm}, we deduce that there exists 
$B_m^{i,j}$ such that $X_m^{i,j} = \nabla^{\perp} B_m^{i,j}$. 
Now, we still denote by $\tilde{u}_m$ and $A_m$ extensions of $\tilde{u}_m$ and $A_m$ on $L^{\infty}\cap H^1(\mathbb{D})$ such that there is a constant $D$ independent of $m$ with
\begin{equation} \label{ineqext} \left\| \nabla \tilde{u}_m \right\|_{L^2(\mathbb{D})} \leq D \left\| \nabla \tilde{u}_m \right\|_{L^2(\mathbb{D} \setminus \mathbb{D}_{\lambda_m})} \hbox{ and } \left\| \nabla A_m \right\|_{L^2(\mathbb{D})} \leq D \left\| \nabla A_m \right\|_{L^2(\mathbb{D} \setminus \mathbb{D}_{\lambda_m})} \hskip.1cm. \end{equation}
For instance, if $f:\mathbb{D}\setminus \mathbb{D}_{\lambda_m}\rightarrow\mathbb{R}$, we take the extension
$$ f(z) = f\left(\frac{z\lambda_m^2}{\left|z\right|^2}\right) \phi\left(\frac{z\lambda_m^2}{\left|z\right|^2}\right) $$
for $z \in \mathbb{D}_{\lambda_m}$, where $\phi \in \mathcal{C}_c^{\infty}(\mathbb{D})$ is a cut-off function such that $\phi = 1$ on $\mathbb{D}_{\lambda_m}$.

We let $D_m$ be such that 
\begin{equation} \label{eqDm} \left\{ \begin{array}{ll} 
\Delta D^i_m =  \left\langle \nabla^{\perp} A_m , \nabla \tilde{u}_m^i \right\rangle & \hbox{ in } \mathbb{D} \\
D_m = 0 & \hbox{ on } \partial \mathbb{D} 
\end{array} \right. \end{equation} 
Since $rot(A_m \nabla\tilde{u}_m - \nabla^{\perp} D_m) =0 $, there is $C_m$ such that $\nabla C_m = A_m \nabla\tilde{u}_m - \nabla^{\perp} D_m$, and $C_m$ satisfies the equation
$$ \Delta C_m^j = \sum_i \left\langle \nabla^{\perp}B_m^{i,j} , \nabla \tilde{u}_m^j \right\rangle $$
on $\mathbb{D} \setminus \mathbb{D}_{\lambda_m}$. We set $C_m = \psi_m + v_m$ where $v_m$ is harmonic and $\psi_m$ satisfies
\begin{equation} \label{eqpsim} \left\{ \begin{array}{ll} 
\Delta \psi_m^j = \sum_i \left\langle \nabla^{\perp}B_m^{i,j} , \nabla \tilde{u}_m^j \right\rangle & \hbox{ in } \mathbb{D}\setminus \mathbb{D}_{\lambda_m} \\
\psi_m = 0 & \hbox{ on } \partial \mathbb{D} \cup \partial \mathbb{D}_{\lambda_m} \hskip.1cm.
\end{array} \right. \end{equation} 
Wente's estimates involving the $L^{2,1}$-norm on the disc for \eqref{eqDm} with the estimates \eqref{ineqext} on the extensions $\tilde{u}_m$ and $A_m$ and $L^{2,1}$-Wente's estimates on the annulus (see \cite{LAU2014}, lemma 2.1) for \eqref{eqpsim} give constants $K_0$ and $K_1$ independent of $m$ such that
\begin{equation} \label{wenteDm} \left\|\nabla D_m\right\|_{L^{2,1}(\mathbb{D})} \leq K_0 \left\|\nabla \tilde{u}_m\right\|_{L^2(\mathbb{D}\setminus \mathbb{D}_{\lambda_m})}^2 \end{equation}
\begin{equation} \label{wentepsim} \left\|\nabla \psi_m \right\|_{L^{2,1}(\mathbb{D}\setminus \mathbb{D}_{\lambda_m})} \leq K_1 \left\|\nabla \tilde{u}_m\right\|_{L^2(\mathbb{D}\setminus \mathbb{D}_{\lambda_m})}^2 \hskip.1cm. \end{equation}
Since $v_m$ is harmonic, we get a Fourier series
$$ v_m = c_m^0 + d_m^0 \ln(r) + \sum_{p\in\mathbb{Z}^{\star}} \left(c_m^p r^{p} + d_m^p r^{-p}\right) e^{ip\theta} $$
and since $\nabla_{\theta} v_m$ has no logarithm part, we use \cite{LAU2014}, lemma A.2, to get a constant $K_2$ independent of $m$ such that
$$ \left\|\nabla_{\theta} v_m \right\|_{L^{2,1}(\mathbb{D}_{\frac{1}{2}}\setminus \mathbb{D}_{2\lambda_m})} \leq K_2 \left\|\nabla v_m \right\|_{L^{2}(\mathbb{D}\setminus \mathbb{D}_{\lambda_m})} \leq K_2 \left\|\nabla C_m \right\|_{L^{2}(\mathbb{D}\setminus \mathbb{D}_{\lambda_m})} $$
since $v_m$ is the harmonic extension of $C_m$ on $\mathbb{D}\setminus \mathbb{D}_{\lambda_m}$. Then, by Young inequalities and since $A_m \leq 1$, 

\begin{equation} \label{angularestvm} \left\|\nabla_{\theta} v_m \right\|_{L^{2,1}(\mathbb{D}_{\frac{1}{2}}\setminus \mathbb{D}_{2\lambda_m})}^2 \leq 2 K_2^2 \left( \left\|\nabla D_m \right\|_{L^{2}(\mathbb{D}\setminus \mathbb{D}_{\lambda_m})}^2 + \left\|\nabla \tilde{u}_m \right\|_{L^{2}(\mathbb{D}\setminus \mathbb{D}_{\lambda_m})}^2 \right) \end{equation}

Now, (\ref{wenteDm}), (\ref{wentepsim}) and
(\ref{angularestvm}) give that

\begin{eqnarray*}
\left\|\nabla_{\theta} \tilde{u}_m \right\|_{L^{2,1}(\mathbb{D}_{\frac{1}{2}}\setminus \mathbb{D}_{2\lambda_m})} &\leq& 4(K_0 + K_1) \left\|\nabla \tilde{u}_m \right\|_{L^2(\mathbb{D}\setminus \mathbb{D}_{\lambda_m})}^2 \\ 
&&+ 4\sqrt{2} K_2 \sqrt{ K_0^2  \left\|\nabla \tilde{u}_m \right\|_{L^2(\mathbb{D}\setminus \mathbb{D}_{\lambda_m})}^4 +  \left\|\nabla \tilde{u}_m \right\|_{L^2(\mathbb{D}\setminus \mathbb{D}_{\lambda_m})}^2 } \hskip.1cm.
\end{eqnarray*}
Looking at this inequality in ${\mathbb A}_{m,R}$ completes the proof of Step 3.

\medskip Gathering Step 2 and Step 3, the duality $L^{2,1} - L^{2,\infty}$ gives that 
\begin{equation} \label{angular} \lim_{R\to +\infty}\lim_{m\to+\infty} \| \nabla_{\theta} u_m \|_{L^2\left({\mathbb A}_{m,R}\right)} = 0 \hskip.1cm.\end{equation}
Since $u_m$ is harmonic with free boundary in ${\mathbb D}^+$, we have the following Pohozaev identity
\begin{equation}\label{Poh}
 \int_{\partial\mathbb{D}_r^+}  \left\vert \nabla_{\theta} u_m\right\vert^2 =\int_{{\partial \mathbb D_r}^+}\left\vert \nabla_r u_m\right\vert^2  \hbox{ for all }0<r<1\hskip.1cm.
\end{equation}
Indeed, let us write with some integration by parts that 
\begin{eqnarray*}
0&=&\int_{{\mathbb D}_r^+} \left(x\left(u_m\right)_x + y \left(u_m\right)_y\right)\cdot \Delta u_m\\
&=& -r\int_{\partial {\mathbb D}_r^+} \left\vert \nabla_r u_m\right\vert^2 + \int_{\left(-1,1\right)\times\left\{0\right\}} \left(x\left(u_m\right)_x + y \left(u_m\right)_y\right)\cdot \left(u_m\right)_y\\
&& + \int_{{\mathbb D}_r^+} \left\vert \nabla u_m\right\vert^2 + \frac{1}{2}  \int_{{\mathbb D}_r^+} \left( x\left(\left\vert \nabla u_m\right\vert^2\right)_x + y \left(\left\vert \nabla u_m\right\vert^2\right)_y\right)\\
&=& -r\int_{\partial {\mathbb D}_r^+} \left\vert \nabla_r u_m\right\vert^2 + \int_{\left(-1,1\right)\times\left\{0\right\}} \left(x\left(u_m\right)_x + y \left(u_m\right)_y\right)\cdot \left(u_m\right)_y+\frac{1}{2}r \int_{\partial {\mathbb D}_r^+} \left\vert \nabla u_m\right\vert^2\hskip.1cm.
\end{eqnarray*}
Now we have that 
$$\int_{\left(-1,1\right)\times\left\{0\right\}} \left(x\left(u_m\right)_x + y \left(u_m\right)_y\right)\cdot \left(u_m\right)_y =0 $$
since on $\left(-1,1\right)\times\left\{0\right\}$, 
$$\left(x\left(u_m\right)_x + y \left(u_m\right)_y\right) \in T_u {\mathbb S}^n$$
and $\left(u_m\right)_y$ is orthogonal to $T_u{\mathbb S}^n$. This proves \eqref{Poh} since $\left\vert \nabla u_m\right\vert^2 = \left\vert \nabla_r u_m\right\vert^2 + \left\vert \nabla_\theta u^m\right\vert^2$. Integrating \eqref{Poh} gives that 
$$\int_{{\mathbb A}_{m,R}} \left\vert \nabla u_m\right\vert^2 = 2 \int_{{\mathbb A}_{m,R}} \left\vert \nabla_\theta u_m\right\vert^2 $$
so that (\ref{limgrad}) follows thanks to (\ref{angular}). 

Finally, we have that
$$ \int_{I_{m,R}}  u_m . (-\partial_t u_m) = \int_{{\mathbb A}_{m,R}} \left|\nabla u_m\right|^2 - \int_{\partial\mathbb{D}_{\frac{1}{R}}^+} u_m . \partial_{r} u_m  + \int_{\partial\mathbb{D}_{R \lambda_m}^+} u_m . \partial_{r} u_m $$
so that
$$\left| \int_{I_{m,R}}  u_m . (-\partial_t u_m) \right| \leq  \int_{{\mathbb A}_{m,R}} \left|\nabla u_m\right|^2 + \int_{\partial\mathbb{D}_{\frac{1}{R}}^+} \left| \nabla u_m \right| + \int_{\partial\mathbb{D}_{R \lambda_m}^+} \left| \nabla u_m \right| $$
and \eqref{limboundary} follows from \eqref{limgrad} and \eqref{limdeltamR}.

This ends the proof of Claim 5, and as already said finishes the proof of Theorem \ref{t2}.

\hfill $\diamondsuit$

\bibliographystyle{plain}
\bibliography{biblio}

\end{document}